\documentclass[12pt,twoside]{article}
\usepackage{amssymb}
\usepackage{amsmath}
\usepackage{amsfonts}
\usepackage{tikz}
\usetikzlibrary{arrows}

\usepackage{amsthm}
\usepackage{calc}

\usepackage{float}
\usepackage{newlfont}
\usetikzlibrary{calc}

\date{}

\begin{document}

\vspace*{3cm}

\centerline{}

\centerline {\Large{\bf Comment on ``Two notes on}}
\centerline {\Large{\bf imbedded prime divisors" }}

\centerline{}

\centerline{\bf {Rahul Kumar\footnote{The author was supported by the SRF grant from UGC India, Sr.
No. 2061440976.} \& Atul
Gaur\footnote{The author was supported by the MATRICS grant from DST-SERB India, No. MTR/2018/000707.}}}

\centerline{Department of Mathematics}

\centerline{University of Delhi, Delhi, India.}

\centerline{E-Mail: rahulkmr977@gmail.com; gaursatul@gmail.com}

\centerline{}

\newtheorem{Theorem}{\quad Theorem}[section]

\newtheorem{Corollary}[Theorem]{\quad Corollary}

\newtheorem{Lemma}[Theorem]{\quad Lemma}

\newtheorem{Proposition}[Theorem]{\quad Proposition}

\theoremstyle{definition}

\newtheorem{Definition}[Theorem]{\quad Definition}

\newtheorem{Example}[Theorem]{\quad Example}

\newtheorem{Remark}[Theorem]{\quad Remark}

\begin{abstract} In this note, we show that a part of \cite[Remark~2.2]{ratliff} is not correct. Some conditions are given under which the same holds.      
\end{abstract}


\vspace{0.5cm}

The following result was proved in $\cite[Remark~2.2]{ratliff}$.

\begin{Theorem}\label{t1}
If $R\subset T$ are Noetherian rings such that there does not exist any integrally dependent adjacent Noetherian rings between them, then for each $\bar{c}/\bar{b}\in T/Z$ (where $Z = Rad(T) = Rad(R)$ and $\bar{b}, \bar{c}$ regular in $R/Z$), we have either $\bar{c}/\bar{b}\in R/Z$ or $\bar{b}/\bar{c}\in R/Z$, and so $(R/Z)[\bar{c}/\bar{b}]$ is a localization of $R/Z$.  
\end{Theorem}

The following examples show that the above result is not correct:

\begin{Example}\label{e1}
\begin{enumerate}
\item[(1)]Let $R = \mathbb{Z}$ and $T = \mathbb{Z}[1/p]$, where $p$ is a prime number. Then clearly $R\subset T$ are adjacent Noetherian rings, $R$ is integrally closed in $T$, and $Rad(R) = Rad(T) =  0$. Now, $q/p\in T\setminus R$ for any prime $q$ distinct from $p$ but $p/q\notin R$.

\item[(2)]Let $R = \mathbb{Z}(+)M$ and $T = \mathbb{Z}[1/p](+)M$, where $p$ is a prime number and $M$ is any finitely generated $R$-module as well as $T$-module. Clearly, $R$ is integrally closed in $T$ and $R\subset T$ is adjacent. Also, by \cite[Theorem~4.8]{ander}, $R\subset T$ are Noetherian rings. Note that $Z = Rad(R) = Rad(T) =  0(+)M$ and so by \cite[Theorem~3.1]{ander}, $R/Z = \mathbb{Z}$ and $T/Z = \mathbb{Z}[1/p]$. Now, $q/p\in T/Z$ for any prime $q$ distinct from $p$ but neither $q/p\notin R/Z$ nor $p/q\notin R/Z$. 
\end{enumerate}
\end{Example}

\begin{Remark}
\begin{enumerate}
\item[(1)] Let $R\subset T$ be an extension of Noetherian rings such that there does not exist any integrally dependent adjacent Noetherian rings between them. Then $R$ is integrally closed in $T$, by the proof of \cite[Theorem~2.1]{ratliff}, and hence $Rad(T) = Rad(R)$. In addition, if $Z(T) = Z(R)$, then $Z(T/Rad(R)) = Z(R/Rad(R))$ as the set of zero-divisors in a reduced ring is the union of minimal prime ideals. Now, if Theorem \ref{t1} holds, then it is easy to see that $R/Rad(R)$ and $T/Rad(R)$ are local (proof follows mutatis mutandis from the proof of \cite[Theorem~1.5]{ayache} and \cite[Corollary~1.6]{ayache}). Consequently, $R$ and $T$ are local. Conversely, if $R\subset T$ is an extension of local Noetherian rings such that there does not exist any integrally dependent adjacent Noetherian rings between them, then the proof of \cite[Theorem~2.1]{ratliff} shows that Theorem \ref{t1} holds. In particular, for an extension of domains, Theorem \ref{t1} holds if and only if the domains are local. 
\item[(2)] Note that if $(R,T)$ is a normal pair (in the sense of \cite{davis}), then there does not exist any integrally dependent adjacent rings between them. Let $R\subset T$ be an extension of Noetherian rings such that $(R,T)$ is a normal pair. Then $T/Rad(R)$ is an overring of $R/Rad(R)$, by the proof of \cite[Theorem~2.1]{ratliff}. Now, if there is at most one maximal ideal of $R$ which is not minimal prime, then the proof of Theorem \ref{t1} follows mutatis mutandis from the proof of \cite[Proposition~3.6(b)]{dobbs2} and part $(1)$. Note that the condition that $R$ is complemented in \cite[Proposition~3.6(b)]{dobbs2} is not required here. 
\end{enumerate}
\end{Remark}


\end{document}